\documentclass[12pt]{article}
\usepackage{amsmath,amssymb,a4wide}

 \def\Gal{\mathop{\rm Gal}\nolimits}
 \def\End{\mathop{\rm End}\nolimits}

 \def\Spec{\mathop{\rm Spec}\nolimits}
 
 \def\deg{\mathop{\rm deg}\nolimits}

\def\ord{\mathop{\rm ord}\nolimits}

\def\GL{\mathop{\rm GL}\nolimits}

\def\sep{{\rm sep}}

\let\phi\varphi
\let\epsilon\varepsilon

\newtheorem{Thm}{Theorem}[section]
\newtheorem{Prop}[Thm]{Proposition}
\newtheorem{Lem}[Thm]{Lemma}

\def\qed{{\hskip0pt\unskip\unskip\nobreak\hfil\penalty50
          \hskip1em\hbox{}\nobreak\hfil
           {$\square$}
          \parfillskip=0pt\finalhyphendemerits=0
          \par}\medskip}

               {\noindent{\bf Proof.}\ }
               {\qed}

               {\noindent{\bf Proof of #1.}\ }
               {\qed}

\newcommand{\BC}{{\mathbb{C}}}

\newcommand{\BF}{{\mathbb{F}}}
\newcommand{\BG}{{\mathbb{G}}}

\newcommand{\BP}{{\mathbb{P}}}
\newcommand{\BQ}{{\mathbb{Q}}}

\newbox\mybox
\def\arrover#1{\mathrel{
       \setbox\mybox=\hbox spread 1.4em
              {\hfil$\scriptstyle#1$\hfil}
       \vbox{\offinterlineskip\copy\mybox
             \hbox to\wd\mybox{\rightarrowfill}}}}

\def\larrover#1{\mathrel{
       \setbox\mybox=\hbox spread 1.4em
              {\hfil$\scriptstyle#1\vphantom{g}$\hfil}
       \vbox{\offinterlineskip\copy\mybox
             \hbox to\wd\mybox{\leftarrowfill}}}}

\def\ontoover#1{\mathrel{
       \setbox\mybox=\hbox spread 1.4em
              {\hfil$\scriptstyle#1\vphantom{g}$\hfil}
       \vbox{\offinterlineskip\copy\mybox
             \hbox to\wd\mybox{\rightarrowfill\hskip-2.8mm
                               $\rightarrow$}}}}
\def\leftontoover#1{\mathrel{
       \setbox\mybox=\hbox spread 1.4em
              {\hfil$\scriptstyle#1\vphantom{g}$\hfil}
       \vbox{\offinterlineskip\copy\mybox
             \hbox to\wd\mybox{$\leftarrow$\hskip-2.8mm
                               \leftarrowfill}}}}

\def\Cinf{{\BC}_\infty}

\begin{document}

\title{Drinfeld Modular Polynomials in Higher Rank}
\author{Florian Breuer \and Hans-Georg R\"uck}
%\date{}
\maketitle

\begin{abstract}
We study modular polynomials classifying cyclic isogenies between Drinfeld modules of arbitrary rank over the ring $\BF_q[T]$.
\end{abstract}

%%%%%%%%%%%%%%%%%%%%%%%%%%%%%%%
%%%%%%%%%%%%%%%%%%%%%%%%%%%%%%%
\section{Introduction and Results}

Let $F$ be an algebraically closed field which contains the ring
$A:=\mathbb F_q[T]$. We denote by $\End_{\BF_q}(\BG_{a,F})$ the ring of $\BF_q$-linear endomorphisms of the additive group over $F$,
it is isomorphic to the non-commutative ring of $\BF_q$-linear polynomials in $X$ with coefficients in $F$ and multiplication defined
by composition of polynomials.

A Drinfeld module $\rho :A \rightarrow \End_{\BF_q}(\BG_{a,F})$ of rank $r$ in generic characteristic is uniquely determined by
\[
\rho(T)(X) = T X + g_1(\rho) X^q + \cdots + g_{r-1}(\rho) X^{q^{r-1}} + \Delta(\rho) X^{q^r}
\]
with coefficients $g_1(\rho), \ldots, g_{r-1}(\rho), \Delta(\rho) \in F$ and $\Delta(\rho) \neq 0$.
We refer the reader to \cite[chapter 4]{GossBS} for an overview of Drinfeld modules.

The coefficients $g_1(\rho), \ldots , g_{r-1}(\rho), \Delta(\rho)$
describe the isomorphism class of $\rho$ in the following way
(see \cite{Potemine}). Let $g_1 , \ldots , g_{r-1}, \Delta$ be indeterminants
and define
\[
j_k =
\frac{g_k^{(q^r-1)/(q^{\gcd(k,r)}-1)}}{\Delta^{(q^k-1)/(q^{\gcd(k,r)}-1)}}
\text{ for } k = 1, \ldots, r-1 .
\]
Let in addition $u_1, \ldots , u_{r-1}$ satisfy the Kummer equations
\[
u_k^{(q^r-1)/(q^{\gcd(k,r)}-1)} = j_k \text{ for } k = 1,
\ldots, r-1.
\]
Consider the ring $A[u_1, \ldots ,u_{r-1}]$. The group $G = \mathbb F_{q^r}^\ast / \mathbb F_q^\ast$
acts on it by
\[
u_k^{\bar{\beta}} = \beta^{q^k-1}\cdot u_k \text{ for } \beta \in \mathbb F_{q^r}^\ast.
\]
The subring of invariant elements, which is denoted by
$A[u_1, \ldots , u_{r-1}]^G$, plays the crucial part in the following
isomorphism problem:\\
Two Drinfeld modules $\rho$ and $\tilde{\rho}$ are isomorphic over
$F$ if and only if $I(\rho) = I(\tilde{\rho})$ for each invariant
$I \in A[u_1, \ldots , u_{r-1}]^G$. In other words
the affine space
\[
M^r = \Spec(A[u_1, \ldots , u_{r-1}]^G)
\]
is the coarse moduli space for Drinfeld modules of
rank $r$ and no level structure.

Let $\Cinf$ denote the completion of an algebraic closure of $\BF_q((\frac{1}{T}))$, and denote by
\[
\Omega^r = \BP^{r-1} \smallsetminus \{\text{Linear subvarieties defined over $\BF_q(({\scriptstyle \frac{1}{T}}))$}\}
\]
the Drinfeld upper half-space. Then the moduli space $M^r$ is given analytically by (see \cite{GekelerDMC})
\[
M^r(\Cinf) \cong \GL_r(A)\backslash\Omega^r.
\]

%({\bf rausschmeissen? :})
%\[
%\xymatrix{
%& \Spec(A[u_1,\ldots,u_{r-1}])\ar@{-}[d]\ar@{=}[r] & \BA_A^{r-1} \\
%M^r\ar@{=}[r] & \Spec(A[u_1,\ldots,u_{r-1}]^G)\ar@{-}[d] & \\
%& \Spec(A[j_1,\ldots,j_{r-1}])\ar@{=}[r] & \BA_A^{r-1}
%}
%\]

Let again $\rho : A \rightarrow \End_{\BF_q}(\BG_{a,F})$ be a Drinfeld module
of rank $r$. For $n \in A$ consider all the isogenies
$\rho \rightarrow \rho_j$, $j \in J(n)$, which are cyclic of degree
$n$, i.e. whose kernel is isomorphic to $A / nA$.
The invariants $I(\rho_j)$ for $I \in A[u_1, \ldots , u_{r-1}]^G$ are algebraic over
$A[u_1, \ldots , u_{r-1}]^G$, and again one gets indeterminants $I_j$ such that
$I_j(\rho) = I(\rho_j)$.

Let $I \in A[u_1, \ldots , u_{r-1}]^G$ be an invariant,
then the polynomial
\[
P_{I,n}(X) := \prod_{j \in J(n)} (X - I_j)
\]
is called the {\em modular polynomial} of the invariant $I$ and the level
$n$.
%It is not difficult to see that $P_{I,n}(X) \in A[u_1, \ldots , u_{r-1}]^G[X]$,
%i.e. that its coefficients are again invariants.

Let $f\in F[u_1,\ldots,u_{r-1}]$, we denote by $w(f)$ the weighted degree of $f$, where each monomial is assigned the weight
\[
w(u_1^{\alpha_1}\cdots u_{r-1}^{\alpha_{r-1}}) := \sum_{k=1}^{r-1}\alpha_k\frac{q^k-1}{q^r-1}.
\]
Let $K$ be the quotient field of $A$, i.e. $K=\BF_q(T)$. For a
polynomial $n \in A$ we denote $|n| = q^{\deg n}$.

The aim of this paper is to prove the following result:

\begin{Thm}\label{Main}
%Let $I \in A[u_1, \ldots , u_{r-1}]^G$
%be an invariant with (weighted ?) degree $d$ in the variables
%$u_1, \ldots , u_{r-1}$. Then the (weighted ?) degrees of the
%coefficients of $P_{I,n}(X)$ are bounded by $...$.
Let $I\in A[u_1,\ldots,u_{r-1}]^G$ be an invariant of weighted degree $w(I)$, and $n\in A$ monic. Then
\begin{enumerate}
\item We have $P_{I,n}(X) \in K[u_1,\ldots,u_{r-1}]^G[X]$, which has degree
\[
\#J(n) = |n|^{r-1}\prod_{p|n}\frac{|p|^r}{|p|^r-|p|^{r-1}}
\]
in $X$, and is irreducible in $\Cinf[u_1,\ldots,u_{r-1}][X]$.
\item The weighted degree of the coefficient $a_i\in K[u_1,\ldots,u_{r-1}]^G$ of $X^i$ in $P_{I,n}(X)$ is bounded by:
\[
w(a_i) \leq \left(|n|^{2(r-1)}\prod_{p|n}\frac{|p|^r}{|p|^r-|p|^{r-1}} -i\right)w(I).
\]
\end{enumerate}
\end{Thm}

%{\bf (Can we replace $K$ by $A$?)}

\paragraph{Example.} When $r=2$ we need only the usual $j$-invariant,
\[
j = u_1^{q+1} = \frac{g_1^{q+1}}{\Delta}.
\]
We have
\[
P_{j,n}(X) = \Phi_n(X,j) = \Phi_n(X,u_1^{q+1}),
\]
where $\Phi_n(X,Y) \in A[X,Y]$ is the modular polynomial of level $n$ constructed by Bae  \cite{Bae}.

\paragraph{Outline of the paper.} 
We define in \S2.2 a parameter $q_{\Lambda_{r-1}}(z_r)$, which may be viewed as a local parameter at the cusp of the compactification of $M^r$, although we make no use of this interpretation. Section 2 is devoted to computing the expansions of various lattice invariants in terms of this parameter.

In Section 3 we study sublattices with cyclic quotients, as these correspond to cyclic isogenies. 
After expanding invariants of sublattices in terms of the parameter, we count the number of such sublattices in \S3.2, prove a non-cancelation result for pole orders in \S3.3 and finally prove Theorem \ref{Main} in \S3.4.

\section{Lattice invariants}

In this section we compute various lattice invariants and express them as Laurent series in a parameter ``at the cusp''.

\subsection{Eisenstein Series}\label{Sec1}

Let $F = \Cinf$. There is a one-to-one correspondence between
Drinfeld modules $\rho$ of rank $s$ and $A$-lattices $\Lambda
\subset \mathbb C_\infty$ of rank $s$ given by the following
relation (see e.g. \cite[Chapter 4]{GossBS}): Let
\begin{equation} \label{eq1}
e_\Lambda(z) = z \prod_{0 \neq \lambda \in \Lambda} \left(1 -
\frac{z}{\lambda}\right) = \sum_{i=0}^\infty e_{q^i}(\Lambda) z^{q^i}
\end{equation}
be the exponential function for $\Lambda$. Then the corresponding
Drinfeld module
%$\rho_\Lambda : A \rightarrow \End_{\BF_q}(\BG_{a,\Cinf})$ with
$\rho_\Lambda$ determined by
\begin{equation} \label{eq1a}
\rho_\Lambda(T)(X) = T X + g_1(\Lambda) X^q + \cdots +
g_{s-1}(\Lambda) X^{q^{s-1}} + \Delta(\Lambda) X^{q^s}
\end{equation}
satisfies the equation
\begin{equation} \label{eq2}
e_\Lambda (T z) = \rho_\Lambda(T) \big(e_\Lambda (z)\big) .
\end{equation}
Equation (\ref{eq2}) and the definitions (\ref{eq1}) and
(\ref{eq1a}) yield for each $k \geq 1$ the equation
\begin{equation} \label{eq3}
(T^{q^k} - T) e_{q^k}(\Lambda) = g_k(\Lambda) + \sum_{j=1}^{k-1}
g_j(\Lambda) e_{q^{k-j}}(\Lambda)^{q^j}.
\end{equation}
Hence there is a polynomial $F_k \in A[X_1, \ldots
X_k]$, which is independent of the lattice $\Lambda$ and which can
be computed recursively, such that
\begin{equation} \label{eq4}
g_k(\Lambda) = F_k\big(e_q(\Lambda), \ldots , e_{q^k}(\Lambda)\big).
\end{equation}

Since $e_\Lambda ^\prime (z) = 1$, we get for the logarithmic
derivative of $e_\Lambda (z)$
\begin{equation} \label{eq5}
\frac{1}{e_\lambda (z)} = \sum_{\lambda \in \Lambda}
\frac{1}{z-\lambda} = \frac{1}{z} - \sum_{i=0}^\infty \left(\sum_{0
\neq \lambda \in \Lambda} \lambda^{-(i+1)}\right) z^i .
\end{equation}
Let
\[
E_k(\Lambda) = \sum_{0 \neq \lambda \in \Lambda}\frac{1}{\lambda^{k}}
\]
be the $k$-th {\em Eisenstein series} of $\Lambda$.
$E_k(\Lambda) = 0$ if $k$ is not divisible by $q-1$. Moreover,
we have $E_{kq}(\Lambda) = E_k(\Lambda)^q.$ Using the expansion
(\ref{eq1}) of $e_\Lambda(z)$ and equation (\ref{eq5}) we get
\begin{equation} \label{eq6}
\big(\sum_{i=0}^\infty e_{q^i}(\Lambda) z^{q^i-1}\big) \big(\sum_{j=0}^\infty
E_j(\Lambda) z^j\big) = -1.
\end{equation}
Here we define $E_0(\Lambda) = -1$. Comparing coefficients yields
in particular for each $k \geq 1$
\begin{equation} \label{eq7}
e_{q^k}(\Lambda) = E_{q^k-1}(\Lambda) + \sum_{i=1}^{k-1}
e_{q^i}(\Lambda) E_{q^{k-i}- 1}(\Lambda)^{q^i}.
\end{equation}
Here we see that there is a polynomial $G_k \in \mathbb F_q[X_1,
\ldots , X_k]$, which is independent of the lattice $\Lambda$ and
which can be computed recursively, such that
\begin{equation} \label{eq8}
e_{q^k}(\Lambda) = G_k\big(E_{q-1}(\Lambda), \ldots ,
E_{q^k-1}(\Lambda)\big).
\end{equation}
Equations (\ref{eq4}) and (\ref{eq8}) show that for each $k \geq
1$ there is a polynomial $H_k \in A[X_1, \ldots X_k]$
which is independent of the lattice $\Lambda$ such that
\begin{equation} \label{eq9}
g_k(\Lambda) = H_k\big(E_{q-1}(\Lambda), \ldots , E_{q^k-1}(\Lambda)\big).
\end{equation}
\paragraph{Remark.} The expression (\ref{eq9}) can also be obtained using recursion
formulas involving $g_i$ and $E_{q^i-1}$ directly (see e.g. \cite[\S II.2]{GekelerDMC}).

\subsection{Parameter at the cusp}

The following calculations are deeply influenced by Goss
\cite{GossES}.

Let $\Lambda_{r-1}$ be an $A$-lattice in $\mathbb C_\infty$ of rank
$r-1$ with $\Delta(\Lambda_{r-1})=1$. In addition let $z_r \in
\mathbb C_\infty$ such that $\Lambda_r = \Lambda_{r-1} \oplus A z_r$
is a lattice of rank $r$. We consider
$$q_{\Lambda_{r-1}}(z_r) = \frac{1}{e_{\Lambda_{r-1}}(z_r)}$$
as a local parameter ``at the cusp''. Our first aim is to expand
various lattice functions of $\Lambda_r$ in terms of this
parameter.

To begin with the calculations, we take the logarithmic derivative
of $e_{\Lambda_{r-1}}(z_r)$ and get
$$q_{\Lambda_{r-1}}(z_r) =
\frac{e^\prime_{\Lambda_{r-1}}(z_r)}{e_{\Lambda_{r-1}}(z_r)} =
\sum_{\lambda \in \Lambda_{r-1}}\frac{1}{\lambda + z_r}.$$ In this
section we want to express for each $i$ the sums $\sum_{\lambda
\in \Lambda_{r-1}}(\lambda + z_r)^{-i}$ in terms of
$q_{\Lambda_{r-1}}(z_r)$. Let $x$ be arbitrary, we get
\begin{equation} \label{eq10}
\frac{1}{e_{\Lambda_{r-1}}(x) - e_{\Lambda_{r-1}}(z_r)} =
\frac{1}{e_{\Lambda_{r-1}}(x-z_r)} = \sum_{\lambda \in
\Lambda_{r-1}}\frac{1}{x - z_r - \lambda}.
\end{equation}
Using the definition of the local parameter this yields
\begin{equation} \label{eq11}
-q_{\Lambda_{r-1}}(z_r) \, \frac{1}{1 - e_{\Lambda_{r-1}}(x) \,
q_{\Lambda_{r-1}}(z_r)} = \sum_{\lambda \in
\Lambda_{r-1}}\frac{1}{x - (z_r + \lambda)}.
\end{equation}
The Taylor expansion in terms of $x$ of the right side of
(\ref{eq11}) equals
$$\sum_{i=0}^\infty\left(-\sum_{\lambda \in \Lambda_{r-1}}\left(\frac{1}{z_r +
\lambda}\right)^{i+1}\right)\, x^i.$$ Let the left side of (\ref{eq11}) be
given with
\begin{equation} \label{eq12}
\frac{1}{1 - e_{\Lambda_{r-1}}(x)\, q_{\Lambda_{r-1}}(z_r)}=
\sum_{i=0}^\infty A_i \, x^i,
\end{equation}
then comparing both sides of (\ref{eq11}) yields
\begin{equation} \label{eq13}
\sum_{\lambda \in \Lambda_{r-1}}\left(\frac{1}{z_r + \lambda}\right)^{i+1} =
q_{\Lambda_{r-1}}(z_r) \, A_i.
\end{equation}
Therefore we have to expand $A_i$ in terms of the local parameter.
Let as before
$$e_{\Lambda_{r-1}}(x) = \sum_{j=0}^\infty e_{q^j}(\Lambda_{r-1})\, x^{q^j}.$$
Then (\ref{eq12}) says
$$\big(\sum_{i=0}^\infty A_i \, x^i\big)\left(1 - q_{\Lambda_{r-1}}(z_r) \sum_{j=0}^\infty
e_{q^j}(\Lambda_{r-1})\, x^{q^j}\right) = 1.$$ This yields  $A_0 = 1$ and
the following recurrences for $i \geq 1$
\begin{equation} \label{eq14}
A_i = q_{\Lambda_{r-1}}(z_r) \, A_{i-1} + q_{\Lambda_{r-1}}(z_r)
\sum_{j \geq 1, q^j \leq i} A_{i-q^j} \, e_{q^j}(\Lambda_{r-1}).
\end{equation}
We see that $A_i = P_i(q_{\Lambda_{r-1}}(z_r))$ where $P_i$ is a
monic polynomial of degree $i$, which is divisible by
$q_{\Lambda_{r-1}}(z_r)$ if $i \geq 1$ and whose coefficients are
elements of $\mathbb F_q[e_{q^j}(\Lambda_{r-1})\mid q^j \leq i]$.
Hence (\ref{eq13}) can be written as
\begin{equation} \label{eq15}
\sum_{\lambda \in \Lambda_{r-1}} \left(\frac{1}{z_r + \lambda}\right)^{i+1} =
q_{\Lambda_{r-1}}(z_r) P_i\big(q_{\Lambda_{r-1}}(z_r)\big).
\end{equation}

\subsection{$q_{\Lambda_{r-1}}(z_r)$-expansions of Eisenstein Series}

In this section we want to expand the Eisenstein series
$E_k(\Lambda_r)$ and the Drinfeld coefficients $g_k(\Lambda_r)$ in
terms of the local parameter $q_{\Lambda_{r-1}}(z_r)$. We
calculate
\begin{eqnarray}
E_k(\Lambda_r) & = & E_k(\Lambda_{r-1}) + \sum_{0 \neq a \in A}
\sum_{\lambda \in \Lambda_{r-1}} \left(\frac{1}{\lambda
+ a z_r}\right)^k \nonumber\\ \label{eq16} & = & E_k(\Lambda_{r-1}) +
\sum_{0 \neq a \in A} q_{\Lambda_{r-1}}(a z_r)
P_{k-1}\big(q_{\Lambda_{r-1}}(a z_r)\big)
\end{eqnarray}
using (\ref{eq15}) with the polynomials $P_{k-1}$.

Let $\rho_{\Lambda_{r-1}}$ be the Drinfeld module corresponding to
the lattice $\Lambda_{r-1}$. Since $\Delta(\Lambda_{r-1}) = 1$ by
assumption, we get for $a \in A$ with leading
coefficient $l(a)$
$$\rho_{\Lambda_{r-1}}(a)(X) = a X + \ldots + l(a) X^{q^{(r-1) \deg a,}}$$
where all the coefficients are elements in
$A[g_1(\Lambda_{r-1}), \ldots , g_{r-2}(\Lambda_{r-1})]$. The
fundamental relation
$$q_{\Lambda_{r-1}}(a z_r)^{-1} = e_{\Lambda_{r-1}}(a z_r) =
\rho_{\Lambda_{r-1}}(a)\big(e_{\Lambda_{r-1}}(z_r)\big) =
\rho_{\Lambda_{r-1}}(a)\big(q_{\Lambda_{r-1}}(z_r)^{-1}\big)$$ yields a
power series expansion
\begin{equation} \label{eq17}
q_{\Lambda_{r-1}}(a z_r) = l(a)^{-1}
q_{\Lambda_{r-1}}(z_r)^{q^{(r-1)\deg a}} + \sum_{i > q^{(r-1)\deg
a}} a_i(\Lambda_{r-1})\, q_{\Lambda_{r-1}}(z_r)^i
\end{equation}
where the coefficients $a_i(\Lambda_{r-1})$ are in
$A[g_1(\Lambda_{r-1}), \ldots , g_{r-2}(\Lambda_{r-1})]$. This
expansion (\ref{eq17}), equation (\ref{eq16}) and the properties
of the polynomials $P_k$ yield the following expansion of the
Eisenstein series
\begin{equation} \label{eq18}
E_k(\Lambda_r) = E_k(\Lambda_{r-1}) + \sum_{i=1}^\infty
b_i^{(k)}(\Lambda_{r-1})\, q_{\Lambda_{r-1}}(z_r)^i
\end{equation}
where the coefficients $b_i^{(k)}(\Lambda_{r-1})$ are elements of
$A[e_q(\Lambda_{r-1}), \ldots ,
e_{q^{r-2}}(\Lambda_{r-1})]$. Here we used in addition the fact,
that the $g_i$ are polynomials in $e_{q^i}$ (cf. (\ref{eq4})).

In \S\ref{Sec1} (cf. (\ref{eq9})) we saw that the elements $g_k(\Lambda)$
are polynomials\\ $H_k\big(E_{q-1}(\Lambda), \ldots ,
E_{q^k-1}(\Lambda)\big)$ in the Eisenstein series, where the $H_k$'s are
independent of the lattice $\Lambda$. Therefore (\ref{eq18})
yields the expansion
\begin{equation} \label{eq19}
g_k(\Lambda_r) = g_k(\Lambda_{r-1}) + \sum_{i=1}^\infty
c_i^{(k)}(\Lambda_{r-1})\, q_{\Lambda_{r-1}}(z_r)^i
\end{equation}
with coefficients $c_i^{(k)}(\Lambda_{r-1}) \in
A[e_q(\Lambda_{r-1}), \ldots ,e_{q^{r-2}}(\Lambda_{r-1})]$.

\subsection{Expansions of $\Delta(\Lambda_r)$}

We want to expand the discriminant $\Delta(\Lambda_r)$. This can be
done using equation (\ref{eq19}) and the facts $\Delta(\Lambda_r) =
g_r(\Lambda_r)$ and $g_r(\Lambda_{r-1}) = 0$. Then one sees
immediately that the $q_{\Lambda_{r-1}}(z_r)$-order of
$\Delta(\Lambda_r)$ is positive. To get the exact value of this
order one has to evaluate the coefficients
$c_i^{(r)}(\Lambda_{r-1})$ in detail. We proceed in a different way,
using a product expansion of $\Delta(\Lambda_r)$ due to Gekeler
\cite{GekelerPE} and Hamahata \cite{Hamahata}.
%
%Gekeler's product expansion (cf. ?) which we adopt for
%arbitrary rank Drinfeld modules.

The key ingredient is the fundamental relation for each lattice
$\Lambda$ :
\begin{equation} \label{eq20}
\rho_\Lambda(T)(X) = \Delta(\Lambda) \cdot \prod_{z \in
T^{-1}\Lambda/\Lambda} (X - e_\Lambda(z)).
\end{equation}
We get immediately for $\Lambda = \Lambda_r$
\begin{equation} \label{eq21}
\Delta(\Lambda_r) = T \cdot \prod_{0 \neq z \in
T^{-1}\Lambda_r/\Lambda_r} \frac{1}{e_{\Lambda_r} (z)}.
\end{equation}
So it is enough to expand $\prod_{0 \neq z \in
T^{-1}\Lambda_r/\Lambda_r} e_{\Lambda_r} (z)$ in terms of the
local parameter. It is not difficult to show (see \cite[Lemma 1]{Hamahata}) that
\begin{equation} \label{eq22}
e_{\Lambda_r}(z) = e_{\Lambda_{r-1}}(z) \prod_{0 \neq a \in A}
\frac{e_{\Lambda_{r-1}}(z+az_r)}{e_{\Lambda_{r-1}}(a z_r)}.
\end{equation}
We decompose $z \in T^{-1}\Lambda_r/\Lambda_r$ as $z = z^\prime +
\frac{\epsilon}{T}z_r$ with $z^\prime \in
T^{-1}\Lambda_{r-1}/\Lambda_{r-1}$ and $\epsilon \in \mathbb F_q$.
In view of (\ref{eq21}) and (\ref{eq22}) we have to evaluate
\begin{eqnarray} \label{eq23}
\prod_{0 \neq z \in T^{-1}\Lambda_r/\Lambda_r}
e_{\Lambda_{r-1}}(z+y) &=& \\ \left(\prod_{0\neq z^\prime \in
T^{-1}\Lambda_{r-1}/\Lambda_{r-1}} e_{\Lambda_{r-1}}(z^\prime
+y)\,\right) & \cdot & \prod_{\epsilon \neq 0}\left(\prod_{z^\prime \in
T^{-1}\Lambda_{r-1}/\Lambda_{r-1}} e_{\Lambda_{r-1}}(z^\prime +
\frac{\epsilon}{T}z_r+y)\right) \nonumber
\end{eqnarray}for $y=0$ and $y=a z_r$.
We use the fundamental relation (\ref{eq20}) now for $\Lambda =
\Lambda_{r-1}$, where we always assume $\Delta(\Lambda_{r-1})=1$,
and get
$$\prod_{0\neq z^\prime \in
T^{-1}\Lambda_{r-1}/\Lambda_{r-1}} e_{\Lambda_{r-1}}(z^\prime +y)
=
\frac{\rho_{\Lambda_{r-1}}(T)\big(e_{\Lambda_{r-1}}(y)\big)}{e_{\Lambda_{r-1}}(y)}
= \frac{e_{\Lambda_{r-1}}(Ty)}{e_{\Lambda_{r-1}}(y)},$$ this
equals $T$ if $y=0$, and
$$ \prod_{z^\prime \in T^{-1}\Lambda_{r-1}/\Lambda_{r-1}}
e_{\Lambda_{r-1}}(z^\prime +\frac{\epsilon}{T}z_r+y) =
\rho_{\Lambda_{r-1}}(T)\big(e_{\Lambda_{r-1}}(\frac{\epsilon}{T}z_r+y)\big)
= e_{\Lambda_{r-1}}(\epsilon z_r+Ty).$$ If we apply these two
formulas to (\ref{eq23}) we get with (\ref{eq21})
\begin{eqnarray}
\frac{1}{\Delta(\Lambda_r)} &=& \big(\prod_{\epsilon \neq 0}
e_{\Lambda_{r-1}}(\epsilon z_r)\big) \cdot \prod_{0 \neq a \in
A} \prod_\epsilon \frac{e_{\Lambda_{r-1}}\big((a
T+\epsilon) z_r\big)}{e_{\Lambda_{r-1}}(a z_r)^{q^{r-1}}} \nonumber \\
\label{eq24} &=& - \,
\frac{1}{q_{\Lambda_{r-1}}(z_r)^{q-1}} \cdot \prod_{0 \neq a \in
A} \prod_\epsilon \frac{q_{\Lambda_{r-1}}(a
z_r)^{q^{r-1}}}{q_{\Lambda_{r-1}}\big((a T+ \epsilon) z_r\big)}.
\end{eqnarray}
We have used the fact that $\prod_{\epsilon\neq 0}\epsilon=-1$.
In equation (\ref{eq17}) we got the expansion of
$q_{\Lambda_{r-1}}(b z_r)$ in terms of the parameter
$q_{\Lambda_{r-1}}(z_r)$ for each $b \in A$. If we
use this, we can evaluate
\begin{equation} \label{eq25}
\frac{q_{\Lambda_{r-1}}(az_r)^{q^{r-1}}}{q_{\Lambda_{r-1}}\big((a
T+\epsilon)z_r\big)} = 1 + \sum_{i=1}^\infty d_i(\Lambda_{r-1}) \,
q_{\Lambda_{r-1}}(z_r)^i
\end{equation}
with $d_i(\Lambda_{r-1}) \in A[e_q(\Lambda_{r-1}),
\ldots ,e_{q^{r-2}}(\Lambda_{r-1})]$. Now (\ref{eq24}) and
(\ref{eq25}) give the final result
\begin{equation} \label{eq26}
\Delta(\Lambda_r) = - \,
q_{\Lambda_{r-1}}(z_r)^{q-1} + \sum_{i=q}^\infty
f_i(\Lambda_{r-1}) \, q_{\Lambda_{r-1}}(z_r)^i
\end{equation}
with coefficients $f_i(\Lambda_{r-1}) \in A[e_q(\Lambda_{r-1}), \ldots ,e_{q^{r-2}}(\Lambda_{r-1})]$.

%We remark that $\prod_{\epsilon \neq 0} \epsilon = -1$ if
%$\text{char}\mathbb F_q \neq 2$ and $1$ otherwise.

\subsection{Expansions of $u_k(\Lambda_r)$}\label{Sec5}

Now we want to expand the invariants of Drinfeld modules
in terms of the local parameter. Let as above
$$u_k^{(q^r-1)/(q^{\gcd(k,r)}-1)} = j_k =
\frac{g_k^{(q^r-1)/(q^{\gcd(k,r)}-1)}}{\Delta^{(q^k-1)/(q^{\gcd(k,r)}-1)}}.
$$
For a lattice $\Lambda$ of rank $r$ we set
$$u_k(\Lambda) = \frac{g_k(\Lambda)}{\Delta^{(q^k-1)/(q^r-1)}}.$$
This definition is not canonical, it could be changed by an
$(q^r-1)/(q^k-1)$-th root of unity. But since we are interested in
the invariants rather than the $u_k$'s themselves, our setting is ultimately
independent of the various choices.

If we combine (\ref{eq19}) and (\ref{eq26}) we get for $k = 1,
\ldots , r-1$
\begin{eqnarray}
u_k(\Lambda_r) = (-1)^{(q^k-1)/(q^r-1)} g_k(\Lambda_{r-1}) \,
q_{\Lambda_{r-1}}(z_r)^{-(q-1)(q^k-1)/(q^r-1)} \nonumber \\
\label{eq27} + \sum_{i>-(q-1)(q^k-1)/(q^r-1)}
h_i^{(k)}(\Lambda_{r-1})\, q_{\Lambda_{r-1}}(z_r)^i
\end{eqnarray}
with coefficients $h_i^{(k)}(\Lambda_{r-1}) \in
A[e_q(\Lambda_{r-1}), \ldots ,e_{q^{r-2}}(\Lambda_{r-1})]$.

\section{Sublattices}

In view of the main theorem we want to expand invariants
of sublattices of $\Lambda_r$ in terms of the local parameter.

\subsection{Invariants of sublattices}

Let $n \in A$ be a monic polynomial and let
$\tilde{\Lambda}_r \subset \Lambda_r$ be a sublattice with cyclic
quotient $\Lambda_r/\tilde{\Lambda}_r \simeq A/nA$.
Then $\tilde{\Lambda}_{r-1} := \tilde{\Lambda}_r
\cap \Lambda_{r-1}$ is a sublattice of $\Lambda_{r-1}$ with
$\Lambda_{r-1}/\tilde{\Lambda}_{r-1} \simeq A/n_2A$
 where $n = n_1 \cdot n_2$ is a decomposition into
monic factors. In addition we find a basis element $w_r = n_1 z_r
+ \lambda$ with $\lambda \in \Lambda_{r-1}$ such that
$\tilde{\Lambda}_r = \tilde{\Lambda}_{r-1} \oplus A w_r$.

If we want to expand the invariants of $\tilde{\Lambda}_r$ as in
\S\ref{Sec5} we have to normalize $\tilde{\Lambda}_{r-1}$ such that its
discriminant equals $1$. Choose $\alpha \in \mathbb C_\infty$ such
that $\Delta(\alpha \tilde{\Lambda}_{r-1}) = 1$ and consider the
lattice $\alpha \tilde{\Lambda}_r = \alpha \tilde{\Lambda}_{r-1}
\oplus A \alpha w_r$. Then we get with (\ref{eq27})
for $k=1, \ldots , r-1$
\begin{eqnarray}
u_k(\tilde{\Lambda}_r) = u_k(\alpha \tilde{\Lambda}_r) =
(-1)^{(q^k-1)/(q^r-1)} g_k(\alpha
\tilde{\Lambda}_{r-1}) \, q_{\alpha \tilde{\Lambda}_{r-1}}(\alpha
w_r)^{-(q-1)(q^k-1)/(q^r-1)} \nonumber \\  \label{eq28} +
\sum_{i>-(q-1)(q^k-1)/(q^r-1)} h_i^{(k)}(\alpha
\tilde{\Lambda}_{r-1})\, q_{\alpha \tilde{\Lambda}_{r-1}}(\alpha
w_r)^i.
\end{eqnarray}
We want to find $\alpha$, express the parameter $q_{\alpha
\tilde{\Lambda}_{r-1}}(\alpha w_r)$ in terms of
$q_{\Lambda_{r-1}}(z_r)$ and compute the coefficients $g_k(\alpha
\tilde{\Lambda}_{r-1})$, $h_i^{(k)}(\alpha \tilde{\Lambda}_{r-1})$
with formulas involving $\Lambda_{r-1}$.

Since $\Lambda_{r-1}/\tilde{\Lambda}_{r-1} \simeq A/n_2A$, we have $n_2 \Lambda_{r-1} \subset
\tilde{\Lambda}_{r-1}$. We consider the polynomial
$$P(X) = n_2 \, X \prod_{0 \neq \lambda \in
n_2^{-1}\tilde{\Lambda}_{r-1}/\Lambda_{r-1}} \left(1 -
\frac{X}{e_{\Lambda_{r-1}}(\lambda)}\right).$$ This is a polynomial of
degree $q^{(r-2) \deg n_2} = |n_2^{r-2}|$, and its coefficients are elements of
$A[q_{\Lambda_{r-1}}(\lambda)\mid \lambda \in
n_2^{-1}\Lambda_{r-1}]$. These values $q_{\Lambda_{r-1}}(\lambda)$
are independent of $q_{\Lambda_{r-1}}(z_r)$.\\
$P(X)$ describes the isogeny corresponding to $n_2 \Lambda_{r-1}
\subset \tilde{\Lambda}_{r-1}$ and as usual (cf. \cite[\S4.7]{GossBS}) we get
\begin{equation} \label{eq29}
P\big(e_{\Lambda_{r-1}}(X)\big) = e_{\tilde{\Lambda}_{r-1}}(n_2 \, X)
\end{equation}
and
\begin{equation} \label{eq30}
P\big(\rho_{\Lambda_{r-1}}(T)(X)\big) =
\rho_{\tilde{\Lambda}_{r-1}}(T)\big(P(X)\big).
\end{equation}
Let $c_{P}$ be the leading coefficient of $P(X)$, then
comparing leading coefficients in (\ref{eq30}) yields
$$\Delta(\tilde{\Lambda}_{r-1}) = c_{P}^{-q^{r-1}+1}.$$
Hence we take $\alpha=c_{P}^{-1}$ and calculate
$$\Delta(\alpha \tilde{\Lambda}_{r-1}) = \alpha^{-q^{r-1}+1}
\Delta(\tilde{\Lambda}_{r-1})= 1.$$ Since $q_{\Lambda}(z) =
e_{\Lambda}(z)^{-1}$, equation (\ref{eq29}) can be written as
$$P\left(\frac{1}{q_{\Lambda_{r-1}}(X)}\right) =
\frac{1}{q_{\tilde{\Lambda}_{r-1}}(n_2\, X)},$$ which yields an
expansion
$$q_{\tilde{\Lambda}_{r-1}}(n_2 \, X) = \alpha \,
q_{\Lambda_{r-1}}(X)^{q^{(r-2)\deg n_2}} + \sum_{i > q^{(r-2)\deg
n_2}} k_i \, q_{\Lambda_{r-1}}(X)^i
$$
where $k_i \in A[q_{\Lambda_{r-1}}(\lambda)\mid
\lambda \in n_2^{-1}\Lambda_{r-1}]$. We apply this formula to the
parameter $q_{\alpha \tilde{\Lambda}_{r-1}}(\alpha w_r)$ and get
\begin{eqnarray}
q_{\alpha \tilde{\Lambda}_{r-1}}(\alpha w_r) = \alpha^{-1}
q_{\tilde{\Lambda}_{r-1}}(w_r) = \alpha^{-1}
q_{\tilde{\Lambda}_{r-1}}(n_2 \frac{w_r}{n_2}) \nonumber \\
\label{eq31}
 = q_{\Lambda_{r-1}}(\frac{w_r}{n_2})^{q^{(r-2)\deg n_2}} +
\sum_{i>q^{(r-2)\deg n_2}} \alpha^{-1} k_i \,
q_{\Lambda_{r-1}}(\frac{w_r}{n_2})^i.
\end{eqnarray}
Since $w_r=n_1 z_r + \lambda$ with $\lambda \in \Lambda_{r-1}$, we
get
\begin{eqnarray}
q_{\Lambda_{r-1}}(\frac{w_r}{n_2}) = \left(q_{\Lambda_{r-1}}(\frac{n_1
z_r}{n_2})^{-1} +
q_{\Lambda_{r-1}}(\frac{\lambda}{n_2})^{-1}\right)^{-1} \nonumber
\\ \label{eq32} = q_{\Lambda_{r-1}}(\frac{n_1 z_r}{n_2}) + \sum_{i>1} l_i \,
q_{\Lambda_{r-1}}(\frac{n_1 z_r}{n_2})^i
\end{eqnarray}
with $l_i \in A[q_{\Lambda_{r-1}}(\lambda)\mid
\lambda \in n_2^{-1}\Lambda_{r-1}]$.

On the other hand we calculate with formula (\ref{eq17})
\begin{eqnarray}
q_{\Lambda_{r-1}}(\frac{n_1 z_r}{n_2}) = q_{\Lambda_{r-1}}(n_1^2
\frac{z_r}{n}) \nonumber \\ \label{eq33} =
q_{\Lambda_{r-1}}(\frac{z_r}{n})^{q^{(r-1)\deg n_1^2}} +
\sum_{i>q^{(r-1)\deg n_1^2}} a_i(\Lambda_{r-1}) \,
q_{\Lambda_{r-1}}(\frac{z_r}{n})^i.
\end{eqnarray}
Now applying (\ref{eq32}) and (\ref{eq33}) to (\ref{eq31}) we get
\begin{equation} \label{eq34}
q_{\alpha \tilde{\Lambda}_{r-1}}(\alpha w_r) =
q_{\Lambda_{r-1}}(\frac{z_r}{n})^{|n_1^{2r-2}n_2^{r-2}|} + \sum_{i>|n_1^{2r-2}n_2^{r-2}|} m_i \,
q_{\Lambda_{r-1}}(\frac{z_r}{n})^i
\end{equation}
with coefficients $m_i \in A[q_{\Lambda_{r-1}}(\lambda)\mid \lambda \in
n_2^{-1}\Lambda_{r-1}]$.

The coefficients $g_k(\alpha \tilde{\Lambda}_{r-1})$ and
$h_i^{(k)}(\alpha \tilde{\Lambda}_{r-1})$ in (\ref{eq28}) are
elements of the ring\\ $A[e_q(\alpha
\tilde{\Lambda}_{r-1}), \ldots , e_{q^{r-2}}(\alpha
\tilde{\Lambda}_{r-1})]$. Equation (\ref{eq29}) shows that this is
a subring of \\
$A[e_q(\Lambda_{r-1}), \ldots
,e_{q^{r-2}}(\Lambda_{r-1})][q_{\Lambda_{r-1}}(\lambda) \mid
\lambda \in n^{-1} \Lambda_{r-1}]$.

This remark and equations (\ref{eq28}) and(\ref{eq34}) show that
$u_k(\tilde{\Lambda}_r)$ can be expanded in a series
\begin{eqnarray}
u_k(\tilde{\Lambda}_r)  =  (-1)^{(q^k-1)/(q^r-1)} g_k(\alpha \tilde{\Lambda}_{r-1}) \,
q_{\Lambda_{r-1}}(\frac{z_r}{n})^{-|n_1^{2r-2}n_2^{r-2}|(q-1)(q^k-1)/(q^r-1)} \nonumber \\ \label{eq35}
 +  \sum_{i>-|n_1^{2r-2}n_2^{r-2}|(q-1)(q^k-1)/(q^r-1)}
 r_i^{(k)}(\alpha \tilde{\Lambda}_{r-1})\,q_{\Lambda_{r-1}}(\frac{z_r}{n})^i
\end{eqnarray}
with coefficients $g_k(\alpha \tilde{\Lambda}_{r-1})$ and $
r_i^{(k)}(\alpha \tilde{\Lambda}_{r-1})$ in the ring\\
$A[e_q(\Lambda_{r-1}), \ldots ,
e_{q^{r-2}}(\Lambda_{r-1})][q_{\Lambda_{r-1}}(\lambda) \mid
\lambda \in n^{-1} \Lambda_{r-1}]$.

\subsection{Counting sublattices}
Let $\Lambda_r$ be an $A$-lattice of rank $r$. For $n
\in A$ we denote by $f(n,r)$ the number of
sublattices $\tilde{\Lambda}_r \subset \Lambda_r$ with
$\Lambda_r/\tilde{\Lambda}_r \simeq  A/nA$. Then we get
\begin{Prop}\label{counting}
We have
\[
f(n,r) = |n|^{r-1} \prod_{p|n} \frac{|p|^r-1}{|p|^r-|p|^{r-1}}
\]
where the product is taken over all monic irreducible divisors $p$
of $n$, and where $|m| = q^{\deg m}$ for each $m \in A$.
\end{Prop}

\paragraph{Proof.} It is obvious that $f(n,r)$ is multiplicative,
i.e. $f(n,r) \cdot f(m,r) = f(n \cdot m,r)$ if $\gcd(n,m)=1$. We use
induction to prove the formula for $n=p^s$. For $s=1$ we have to
count the number of upper triangular matrices $M$ with coefficients in
$A$ and determinant $p$ which are in reduced form (this means that the only non-zero entries $x$ above
the diagonal lie above the $p$ and satisfy $|x|<|p|$).  We get
$$ f(p,r) = \sum_{i=1}^r |p|^{i-1} = |p|^{r-1}
\frac{|p|^r-1}{|p|^r-|p|^{r-1}}.$$ Let $s\geq2$ and let
$\tilde{\Lambda}_r \subset \Lambda_r$ with
$\Lambda_r/\tilde{\Lambda}_r \simeq  A/p^sA$. Then there is a unique
lattice $\Lambda'_r$ with $\tilde{\Lambda}_r \subset \Lambda'_r
\subset \Lambda_r$ and $\Lambda_r/\Lambda'_r \simeq  A/p^{s-1}A$. On
the other hand if $\Lambda'_r$ is any sublattice of $\Lambda_r$ with
$\Lambda_r/\Lambda'_r \simeq  A/p^{s-1}A$, we can choose bases of
$\Lambda_r$ and $\Lambda'_r$ such that their connecting matrix is
diagonal of the form $\text{diag}(1,\ldots ,1,p^{s-1})$. Any upper
triangular matrix $M_p$ of determinant $p$ in reduced form gives a
sublattice $\tilde{\Lambda}_r$ of $\Lambda'_r$ of index $p$ with
connecting matrix $M_p \cdot \text{diag}(1,\ldots ,1,p^{s-1})$ to
the chosen basis of $\Lambda_r$. One sees immediately that
$\Lambda_r/\tilde{\Lambda}_r \simeq  A/p^sA$ if and only if the
$(r,r)$-th coefficient of $M_p$ is equal to $p$. Hence we have
$|p|^{r-1}$ choices for $M_p$. Thus we get by induction
$$ f(p^s,r) = |p|^{r-1} f(p^{s-1},r) = |p|^{r-1} |p^{s-1}|^{r-1}
\frac{|p|^r-1}{|p|^r-|p|^{r-1}}$$ which proves the formula for
$n=p^s$.\qed

\subsection{Non-cancelation}

The coefficients of the modular polynomial $P_{I,n}(X)$ are
polynomials in the basic invariants. We want to study the weighted
degree of these polynomials by comparing it with its order as a
series expansion in the local parameter. Therefore we need to show
that the leading terms in these series do not cancel. For this we
will need the following technical lemma.

We consider the homothety classes of lattices $\Lambda_s$ of rank
$s$ as points in the moduli space $M^s(\Cinf)$ equipped with the
analytic topology, i.e. where closed sets are the zero-loci of sets
of analytic functions.

For $k=1, \ldots , r-1$ let $v_k$ be functions on $M^r(\Cinf)$;
suppose that their $q_{\Lambda_{r-1}}(\frac{z_r}{n})$-expansions are
of the form
$$
v_k(\Lambda_{r-1}+A z_r) =
a_k(\Lambda_{r-1})q_{\Lambda_{r-1}}(\frac{z_r}{n})^{-c(q^k-1)/(q^r-1)}
+ \text{higher terms},
$$
where the $a_k$'s are algebraically independent analytic functions
on $M^{r-1}(\Cinf)$ and where $c$ does not depend on $\Lambda_{r-1}$
or on $k$.

Let $f\in\Cinf[X_1,\ldots,X_{r-1}]$ be a polynomial of weighted
degree $w(f)$, where the weighted degree of a monomial is given by
$w(X_1^{\alpha_1} \cdots X_{r-1}^{\alpha_{r-1}}) = \sum_{k=1}^{r-1}
\alpha_k \frac{q^k-1}{q^r-1}$.

\begin{Lem}\label{NC}
There exists a non-empty open subset $S\subset M^{r-1}(\Cinf)$ such
that for any $\Lambda_{r-1} \in S$ the weighted degree and the order
of the $q_{\Lambda_{r-1}}(\frac{z_r}{n})$-expansion satisfy
\[
\ord_{q_{\Lambda_{r-1}}(\frac{z_r}{n})}\big(f(v_1(\Lambda_{r-1}+A
z_r),\ldots,v_{r-1}(\Lambda_{r-1}+A z_r))\big) = -c\; w(f).
\]
\end{Lem}

\paragraph{Proof.}
Let
\[
f(X_1,\ldots,X_{r-1}) =
\sum_{(\alpha_1,\ldots,\alpha_{r-1})}b_{(\alpha_1,\ldots,\alpha_{r-1})}X_1^{\alpha_1}\cdots
X_{r-1}^{\alpha_{r-1}}.
\]
Then the leading term of the
$q_{\Lambda_{r-1}}(\frac{z_r}{n})$-expansion of
$f(v_1(\Lambda_{r-1}+A z_r),\ldots,v_{r-1}(\Lambda_{r-1}+A z_r))$ is
given by
\[
\left(
\sum_{(\alpha_1,\ldots,\alpha_{r-1})}b_{(\alpha_1,\ldots,\alpha_{r-1})}
a_1(\Lambda_{r-1})^{\alpha_1}\cdots
a_{r-1}(\Lambda_{r-1})^{\alpha_{r-1}}\right)
q_{\Lambda_{r-1}}(\frac{z_r}{n})^{-c w(f)},
\]
where the sum is taken over all indices $(\alpha_1,\ldots,\alpha_{r-1})$ satisfying
\[
\sum_{k=1}^{r-1} \alpha_k\frac{q^k-1}{q^r-1} = w(f).
\]
The coefficient is zero only if
$a_1(\Lambda_{r-1}),\ldots,a_{r-1}(\Lambda_{r-1})$ satisfy a
polynomial relation. The locus of $\Lambda_{r-1}\in M^{r-1}(\Cinf)$
for which this relation holds is a proper closed set, since the
$a_k$'s are algebraically independent. The result follows. \qed

{\bf Remark.} We can apply Lemma \ref{NC} to $v_k =u_k$ resp. $v_k =
u_k(\tilde{\Lambda}_{r-1})$ in view of (\ref{eq27}) resp.
(\ref{eq35}) and due to the fact that $g_1,\ldots,g_{r-1}$ are
algebraically independent on $M^{r-1}(\Cinf)$.

\subsection{Proof of the Main Result}

%We start by pointing out that, if $m,n\in A$ with $m|n$, then any function $f$ on $\Omega^r$ which has a
%$q_{\Lambda_{r-1}}(\frac{z_r}{m})$-expansion also has a $q_{\Lambda_{r-1}}(\frac{z_r}{n})$-expansion:
%From
%\begin{eqnarray*}
%e_{\Lambda_{r-1}}(\frac{z_r}{m}) & = & e_{\Lambda_{r-1}}(\frac{n}{m}\frac{z_r}{n}) =
%\rho_{\Lambda_{r-1}}(\frac{n}{m})\big(e_{\Lambda_{r-1}}(\frac{z_r}{n}) \big) \\
%& = & e_{\Lambda_{r-1}}(\frac{z_r}{n})^{q^{(r-1)\deg(n/m)}} + \text{lower terms}
%\end{eqnarray*}
%follows that
%\begin{equation}\label{eq37}
%q_{\Lambda_{r-1}}(\frac{z_r}{m}) = q_{\Lambda_{r-1}}(\frac{z_r}{n})^{|n/m|^{r-1}} + \text{lower terms.}
%\end{equation}

\paragraph{Proof of Theorem \ref{Main}.}
(1) The value for the degree follows from Proposition
\ref{counting}.
%We view the coefficients of $P_{I,n}(X)$ as
%functions on $\Omega^{r}$. They are fixed by $\GL_r(A)$, since this
%group permutes the $I_j$'s. Hence they are functions on $M^r(\Cinf)$
%and thus elements of $\Cinf[u_1,\ldots,u_{r-1}]^G$. {\bf It remains
%to show that the coefficients are in} $A[u_1,\ldots,u_{r-1}]^G$.
Let $\rho$ be a Drinfeld module with invariant $I$. The $I_j$'s correspond to cyclic submodules of
$\rho[n]\cong\big(A/nA\big)^r$ of order $n$, which are permuted transitively by $\GL_r(A)$. It follows that
the coefficients of $P_{I,n}(X)$ are functions on $\GL_r(A)\backslash\Omega^r\cong M^r(\Cinf)$, and so $P_{I,n}(X)$ is an
irreducible polynomial with coefficients in $\Cinf[u_1,\ldots,u_{r-1}]^G$.

We next show how to replace $\Cinf$ by $K$. Let $u_1',\ldots,u_{r-1}'\in K$ be arbitrary. These values correspond to a Drinfeld
module $\rho$ defined over $k$ with $u_k = u_k'$ and invariant $I\in K$. The absolute Galois group $\Gal(K^\sep/K)$ permutes
the set of cyclic submodules of $\rho[n]$ of order $n$, hence permutes the $I_j$'s. Thus the coefficients of $P_{I,n}(X)$,
when specialized to $(u_1',\ldots,u_{r-1}')\in K^{r-1}$, lie in $K$. Since $(u_1',\ldots,u_{r-1}')\in K^{r-1}$ is arbitrary,
it follows that the coefficients of $P_{I,n}(X)$ lie in $K[u_1,\ldots,u_{r-1}]^G$.

(2) Let $\tilde{\Lambda}_r\subset\Lambda_r$ be the sublattice
corresponding to $j\in J(n)$, then $I_j = I(\tilde{\Lambda}_r)$.
Equation (\ref{eq35}) and Lemma \ref{NC} show, for suitably chosen
$\Lambda_{r-1}$, that
\[
\ord_{q_{\Lambda_r}(\frac{z_r}{n})} \big(u_k(\tilde{\Lambda}_r)\big)
= -|n_1^{2r-2} n_2^{r-2}|(q-1)\frac{q^k-1}{q^r-1}
\]
(for $k=1,\ldots,r-1$) implies
\[
\ord_{q_{\Lambda_{r-1}}(\frac{z_r}{n})} (I_j) = -(q-1) |n_1^{2r-2}
n_2^{r-2}| w(I)\geq -(q-1)|n|^{2(r-1)}w(I).
\]
Now let $a_i\in A[u_1,\ldots,u_{r-1}]^G$ be the coefficient of $X^i$ in $P_{I,n}(X)$. Then
\[
a_i = (-1)^d\sum_{(j_1,\ldots,j_d)} I_{j_1}\cdots I_{j_d},\qquad\text{where $d=\#J(n)-i$.}
\]
Hence we get
\begin{eqnarray}
\ord_{q_{\Lambda_{r-1}}(\frac{z_r}{n})}(a_i) & \geq &
-(q-1)\big(\#J(n)-i\big)|n|^{2(r-1)}w(I). \label{eq38}
\end{eqnarray}
On the other hand, since $a_i\in A[u_1,\ldots,u_{r-1}]^G$, using
Lemma \ref{NC}, (\ref{eq17}) and (\ref{eq27}), we get
\begin{equation} \label{eq39}
\ord_{q_{\Lambda_{r-1}}(\frac{z_r}{n})}(a_i) = |n|^{r-1}
\ord_{q_{\Lambda_{r-1}}(z_r)}(a_i) = -(q-1)|n|^{r-1}w(a_i).
\end{equation}
Now (\ref{eq38}), (\ref{eq39}) and Proposition \ref{counting} yield
\[
w(a_i) \leq |n|^{2(r-1)} \left(
\prod_{p|n}\frac{|p|^r}{|p|^r-|p|^{r-1}}-i\right)w(I).
\]
\qed

\begin{center}
\rule{8cm}{0.01cm}
\end{center}

\begin{minipage}[t]{8cm}{\small
Department of Mathematical Sciences \\
University of Stellenbosch \\
Stellenbosch, 7600 \\
South Africa \\
fbreuer@sun.ac.za}
\end{minipage}
%
%\hspace{2cm}
\begin{minipage}[t]{8cm}{\small
Fachbereich f\"ur Mathematik und Informatik  \\
Universit\"at Kassel,\\
Kassel, 34132 \\
Germany \\
rueck@mathematik.uni-kassel.de}
\end{minipage}

\end{document}